\documentclass[letterpaper]{article}

\date{\today}

\newcommand{\myauthor}{Benjamin Antieau}
\newcommand{\mytitle}{Cohomological obstruction theory for Brauer classes and the period-index problem}

\title{\mytitle\footnote{This material is based upon work supported by the NSF under Grant No. DMS-0901373.}}
\author{\myauthor}

\usepackage[pdfstartview=FitH,
            pdfauthor={\myauthor},
            pdftitle={\mytitle},
            colorlinks,
            linkcolor=reference,
            citecolor=citation,
            urlcolor=e-mail,
            backref]{hyperref}

\usepackage{rnote}
\usepackage{theorems}
\usepackage{times}
\usepackage[all,cmtip]{xy}
\usepackage{diagxy}
\usepackage{tikz}
\usetikzlibrary{matrix,arrows}

\begin{document}
\maketitle
\begin{abstract}
  \noindent
  Let $U$ be a connected scheme of finite \'etale cohomological dimension in which every finite set of points is contained in an affine open subscheme.
  Suppose that $\alpha$ is a class in $\Hoh^2(U_{\et},\Gm)_{\tors}$. For each positive integer $m$, the $K$-theory
  of $\alpha$-twisted sheaves is used to identify obstructions to $\alpha$
  being representable by an Azumaya algebra of rank $m^2$. The \'etale index of $\alpha$, denoted $eti(\alpha)$, is the least positive integer such that all the
  obstructions vanish. Let $per(\alpha)$ be the order of $\alpha$ in $\Hoh^2(U_{\et},\Gm)_{\tors}$.
  Methods from stable homotopy theory give an upper bound on the \'etale index that depends on the period of $\alpha$ and the \'{e}tale cohomological dimension of $U$; this bound is expressed
  in terms of the exponents of the stable homotopy groups
  of spheres and the exponents of the stable homotopy groups of $B\left(\ZZ/(per(\alpha))\right)$. As a corollary, if $U$
  is the spectrum of a field of finite cohomological dimension $d$, then $eti(\alpha)|per(\alpha)^{\lfloor\frac{d}{2}\rfloor}$, where $\lfloor\frac{d}{2}\rfloor$ is the integer part
  of $\frac{d}{2}$, whenever $per(\alpha)$ is divided neither by the characteristic of $k$ nor by any primes that are small relative to $d$.

\paragraph{Key Words}
Brauer groups, twisted sheaves, higher algebraic $K$-theory, stable homotopy theory.

\paragraph{Mathematics Subject Classification 2000}
Primary: \href{http://www.ams.org/mathscinet/msc/msc2010.html?t=14Fxx&btn=Current}{14F22}, \href{http://www.ams.org/mathscinet/msc/msc2010.html?t=16Kxx&btn=Current}{16K50}.
Secondary: \href{http://www.ams.org/mathscinet/msc/msc2010.html?t=19Dxx&btn=Current}{19D23}, \href{http://www.ams.org/mathscinet/msc/msc2010.html?t=55Qxx&btn=Current}{55Q10},
\href{http://www.ams.org/mathscinet/msc/msc2010.html?t=55Qxx&btn=Current}{55Q45}.
\end{abstract}

\section{Introduction}

\begin{hypothesis}
    Throughout this paper, $U$ denotes a connected scheme of finite \'etale cohomological dimension $d$ having the property that every finite set of points of $U$ is contained in an
    affine open subscheme. For instance, any quasi-projective scheme over a noetherian base satisfies this hypothesis.
\end{hypothesis}

\begin{definition}[see Definition~\ref{def:eti}]
    For $\alpha\in\Hoh^2(U_{\et},\Gm)$, define $eti(\alpha)$ to be the positive generator of the rank map $\K_0^{\alpha,\et}(U)\rightarrow\ZZ$, where $\K^{\alpha,\et}$ denotes
    $\alpha$-twisted \'etale $K$-theory defined in Definition~\ref{def:etalesheaf}.
\end{definition}

This paper is dedicated to proving the following theorem, with the exception of the \textbf{divisibility} property, which is proven in \cite{antieau_cech_2010}.

\begin{theorem}
    Let $\alpha\in\Hoh^2(U_{\et},\Gm)_{\tors}$. Then, $eti(\alpha)$ has the following properties:
    \begin{enumerate}
        \item   \textbf{computability}: in the descent spectral sequence
            \begin{equation*}
                \Eoh_{2}^{s,t}=\Hoh^s(U_{\et},\mathcal{K}_{t}^{\alpha})\Rightarrow\K_{t-s}^{\alpha,\et}(U)
            \end{equation*}
            for $\alpha$-twisted \'etale $K$-theory,
            the integer $eti(\alpha)\in\ZZ\iso\Hoh^0(U_{\et},\mathcal{K}_0^{\alpha})$ is the smallest positive integer such that $d_k^{\alpha}(eti(\alpha))=0$ for all $k\geq 2$,
            where $d_k^{\alpha}$ is the $k$th differential in the spectral sequence;
        \item   \textbf{divisibility}: $per(\alpha)|eti(\alpha)$, where $per(\alpha)$ is the order of $\alpha$ in $\Hoh^2(U_{\et},\Gm)_{\tors}$;
        \item   \textbf{obstruction}: if $\mathcal{A}$ is an Azumaya algebra in the class of $\alpha$, then $eti(\alpha)|deg(\mathcal{A})$, where $deg(\mathcal{A})$, the degree of $\mathcal{A}$, is the positive
            square-root of the rank of $\mathcal{A}$;
        \item   \textbf{bound}: if $per(\alpha)$ is prime to the characteristics of the residue fields of $U$, then
            \begin{equation*}
              eti(\alpha)|\prod_{j\in\{1,\ldots,d-1\}} l_j^\alpha.
            \end{equation*}
            where $l_j^{\alpha}$ is the least common multiple of the exponents of $\pi_j^s$ and $\pi_j^s(B\ZZ/(per(\alpha)))$.
    \end{enumerate}
    In particular, $eti(\alpha)$ is finite even if $\alpha$ is not representable by an Azumaya algebra.
\end{theorem}

The first property is shown in Lemma~\ref{lem:computability}.
The \textbf{obstruction} property is proven in Theorem~\ref{lem:obstruction}, and the \textbf{bound} property is established in Theorem~\ref{thm:etibound}.
An analysis of the integers $l_j^{\alpha}$, together with the \textbf{divisibility} and \textbf{bound}
properties above and the fact that the period and index have the same prime divisors for Brauer classes on a field, gives the following, Theorem~\ref{thm:periodindex}:

\begin{theorem}[Period-\'Etale Index Theorem]\label{thm:intro2}
    Let $k$ be a field, and let $\alpha\in\Hoh^2(k,\Gm)$. Let $S$ be the set of prime divisors of $per(\alpha)$, and suppose that $d=cd_S k<2\min_{q\in S}(q)$. Then,
    \begin{equation*}
        eti(\alpha)|(per(\alpha))^{\lfloor\frac{d}{2}\rfloor},
    \end{equation*}
    where $\lfloor\frac{d}{2}\rfloor$ is the integer part of $\frac{d}{2}$.
\end{theorem}

The theorem should be viewed as a topological version of the period-index conjecture, attributed to \colliot.

\begin{conjecture}[Period-Index Conjecture]\label{conj:periodindexconjecture}
  If $k$ is a field of dimension $d$, then
  \begin{equation*}
    ind(\alpha)|(per(\alpha))^{d-1}
  \end{equation*}
  for all $\alpha\in\Br(k)$, where $ind(\alpha)$ is the square-root of the rank of the unique division algebra representing $\alpha$.
\end{conjecture}

In the conjecture, the dimension might mean either that $k$ is $C_d$, that $k$ is the
function field of a $d$-dimensional algebraic variety over an algebraically closed field, that $k$ is the function field of a $(d-1)$-dimensional variety over a finite field, or that $k$ is
the function field of a $(d-2)$-dimensional variety over a local field. It is not known what the precise statement should be.

However, the conjecture is known to be false if dimension is taken to be the cohomological dimension of the field. For prime powers $l^e$ and $l^f$, with $e\leq f$, a construction of Merkurjev
\cite{merkurjev_kaplansky_1989} can be used to construct a field $k$ with $cd_l(k)=2$, and a division algebra $D$ over $k$ with $per(D)=l^e$ and $ind(D)=l^f$.

For general background on the conjecture and its importance, see \cite{lieblich_period_2007}. It is known to be true in the following cases, where in fact the period and index coincide:
\begin{itemize}
  \item $p$-adic fields, by class field theory;
  \item number fields, by the Brauer-Hasse-Noether theorem;
  \item $C_2$-fields, when $per(\alpha)=2^a 3^b$, by Artin and Harris \cite{artin_brauer_1982};
  \item function fields $k(X)$ of algebraic surfaces $X$ over an algebraically closed field $k$, by de Jong \cite{de_jong_period_2004};
  \item quotient fields $K$ of excellent henselian two-dimensional local domains with residue field $k$ separably closed when $\alpha$ is a class of period prime to the
    characteristic of $k$, by Colliot-Th\'el\`ene, Ojanguren, and Parimala \cite{colliot_quadratic_2002};
  \item fields $l((t))$ of transcendence degree $1$ over $l$, a characteristic zero field of cohomological dimension $1$, by \colliot, P. Gille, and Parimala
    \cite{colliot_arithmetic_2004}.
\end{itemize}

The conjecture is also known in the following situations. Saltman \cite{saltman_division_1997} showed that
\begin{equation*}
    ind(D)|per(D)^2
\end{equation*}
holds for division algebras over the function fields of curves over $p$-adic fields. Lieblich, in \cite{lieblich_period_2009} has shown that this is also true
for the function fields of surfaces over finite fields. Finally, Lieblich and Krashen have established in \cite{lieblich_period_2007} the sharp relation
\begin{equation*}
    ind(D)|per(D)^d
\end{equation*}
for the function fields of curves over $d$-local fields, such as $k((t_1))\cdots((t_d))$, where $k$ is algebraically closed. Moreover, in these examples, the exponent is the best possible.

\paragraph{Acknowledgments}
This paper is part of my Ph.D. thesis,
and I thank Henri Gillet, my thesis advisor at UIC, as well as David Gepner, Christian Haesemeyer, and Brooke Shipley for discussions and support.
Also, the referee of the paper has made a great number of suggestions resulting in the improvement of the exposition.

\section{Stacks of twisted sheaves}

\begin{proposition}[Artin \cite{artin_joins_1971}]
    If $U$ is a scheme such that every finite set of points is contained in some affine open subscheme, then the sheaf cohomology group $\Hoh^2(U_{\et},\Gm)$ is computable by covers (instead of hypercovers);
    that is, $\check{\Hoh}^2(U_{\et},\Gm)\riso\Hoh^2(U_{\et},\Gm)$.
\end{proposition}

\begin{remark}
    The proposition ensures that no information is lost by using only covers in the constructions and theorems below. However, at the expense of another level of detail, all of the material in this paper
    can be modified to apply to any connected scheme of finite cohomological dimension, provided that one uses $1$-hypercovers instead of covers.
    Indeed, for any scheme $U$, the small \'etale site $U_{\et}$ has fiber products and finite products. Therefore, by \cite[Theorem~V.7.4.1]{sga4.2}, $\Hoh^2(U_{\et},\Gm)$ is computable by $1$-hypercovers.
\end{remark}


Two \'etale stacks play a fundamental r\^ole in this paper.
For background on stacks see \cite[Chapter~4]{fga_explained}. The first stack is $\StProj$, the stack of locally free finite
rank coherent modules and isomorphisms. Thus an object in the category of sections $\StProj_V$ on an \'etale map $V\rightarrow U$ is a locally free finite rank coherent $\mathcal{O}_V$-module.
For brevity, such an object will be called a lffr sheaf. Fix a positive integer $n$.
The second stack is the stack $\nSets$ of sheaves of finite and faithful $\mu_n$-sets, where $\mu_n$ is the sheaf of $n$th roots of unity in the \'etale topology. The category of sections $\nSets_V$
consist of sheaves $F$ with a faithful action of $\mu_n|_V$ such that $F$ decomposes into finitely many orbits. Objects will be called $\mu_n$-sheaves. The morphisms are isomorphisms.
Every object of $\nSets_V$ is a disjoint sum of $\mu_n$-torsors. There is a map of stacks, the unit morphism, $i:\nSets\rightarrow\StProj$ obtained by sending a $\mu_n$-torsor to the associated $\Gm$-torsor,
and then taking the sheaf of sections. Disjoint sums are taken to direct sums.

\begin{definition}\label{def:twisting}
    Let $\alpha\in\Hoh^2(U_{\et},\Gm)$, and suppose that $\mathcal{U}=(U_i)_{i\in I}$ is an \'etale cover such that $\alpha$ comes from the \v{C}ech cocycle $(\alpha_{ijk})$,
    where each $\alpha_{ijk}\in\Gamma(U_{ijk},\Gm)$.
    An $\alpha$-twisted coherent $\mathcal{O}_U$-module consists of a coherent $\mathcal{O}_{U_i}$-module $\mathcal{F}_i$ for each $i\in I$, together with isomorphisms
    $\theta_{ij}:\mathcal{F}_i|_{U_{ij}}\riso\mathcal{F}_j|_{U_{ij}}$ such that $\theta_{ki}\circ\theta_{jk}\circ\theta_{ij}=\alpha_{ijk}\in\Gm(U_{ijk})$. For properties of $\alpha$-twisted sheaves,
    see \cite{lieblich_twisted_2008} or \cite{caldararu_derived_2000}.

\end{definition}

The locally free and finite rank $\alpha$-twisted coherent $\mathcal{O}_U$-modules naturally give rise to a stack $\StProjA$, where the sections over $V\rightarrow U$
are the $\alpha|V$-twisted lffr sheaves.

\begin{lemma}\label{lem:representation}
    Let $\alpha\in\Hoh^2(U_{\et},\Gm)$. If $V\rightarrow U$ is \'etale and $V$ is connected, then there is an Azumaya algebra of rank $n^2$ representing $\alpha|V$
    if and only if there is a $\alpha$-twisted lffr sheaf of rank $n$ in $\StProjA_V$.
    \begin{proof}
        See \cite[Proposition~3.1.2.1]{lieblich_twisted_2008}.
    \end{proof}
\end{lemma}

\begin{lemma}\label{lem:linebundles}
    Let $\alpha\in\Hoh^2(U_{\et},\Gm)$, and let $V\rightarrow U$ be an \'etale map. If $\alpha|_V$ is trivial, there is an $\alpha$-twisted lffr rank $1$ sheaf in $\StProjA_V$.
    \begin{proof}
        This follows from \cite[Proposition~3.1.2.1.iv]{lieblich_twisted_2008}.
    \end{proof}
\end{lemma}

Similarly, if $\beta\in\Hoh^2(U_{\et},\mu_n)$, then there is a twisted form $\nSetsB$ of $\nSets$ constructed in the same way as $\StProjA$ is in Definition~\ref{def:twisting}.

\begin{lemma}\label{lem:twistedunit}
    If $\Hoh^2(U_{\et},\mu_n)\rightarrow\Hoh^2(U_{\et},\Gm)$ sends $\beta$ to $\alpha$,
    then the unit map $i:\nSets\rightarrow\StProj$ twists to give a twisted unit map $i^{\beta}:\nSetsB\rightarrow\StProjA$.
    \begin{proof}
        Suppose for simplicity that $\beta$ is defined on the cover $\mathcal{U}=(U_i)_{i\in I}$ by $\beta_{ijk}\in\mu_n(U_{ijk})$. If $F$ is a $\beta$-twisted $\mu_n$-set, then $F_i=F|_{U_i}$ is a $\mu_n$-set
        for all $i\in I$, and there are isomorphisms $\theta_{ij}:F_i\riso F_j$.
        Thus, $i(F_i)$ is a lffr sheaf, and $i(\theta_{ij})$ give isomorphisms $i(F_i)\riso i(F_j)$ such that
        \begin{equation*}
            i(\theta_{ki})\circ i(\theta_{jk})\circ i(\theta_{ij})=\beta_{ijk},
        \end{equation*}
        where now $\beta_{ijk}$ is viewed as a $2$-cocycle in $\Gm$, which is by hypothesis cohomologous to $\alpha$. Thus $i(F_i)$ and $i(\theta_{ij})$ give the data of an $\alpha$-twisted lffr sheaf.
        The details are left to the reader.
    \end{proof}
\end{lemma}

Both stacks $\nSetsB$ and $\StProjA$ are stacks of symmetric monoidal categories in the following sense. Each category of sections is a symmetric monoidal category, under disjoint union and direct
sum respectively, and the restriction is compatible with this structure.

\section{K-theory}\label{sec:k-theory}

\begin{definition}
    There is a functor
    \begin{equation*}
        \K:SymMon\rightarrow \Spt,
    \end{equation*}
    from the category of symmetric monoidal categories and lax functors to spectra. For details, see \cite[Section~1.6]{thomason_symmetric_1995}.
    This $K$-theory is always connective. If $T$ is a symmetric monoidal category, let $\K_n(T)=\pi_n(\K(T))$ for $n\in\ZZ$.
\end{definition}

\begin{example}
    If $R$ is a commutative ring, and if $\StProj_R$ is the symmetric monoidal category of finitely generated projective $R$-modules and isomorphisms, with direct sum,
    then $\K(\StProj_R)$ agrees with Quillen's higher algebraic $K$-theory of $R$ \cite{grayson_higher_1976}. In particular, $\K_0(R)$
    is the usual Grothendieck group of $R$. Similarly, if $X$ is a scheme, and $\StProj_X$ is the category of locally free and finite rank $\mathcal{O}_X$-modules.
    Then the Quillen $Q$-construction $Q\StProj_X$
    of $\StProj_X$ has a natural structure of symmetric monoidal category under direct sum. Quillen's higher algebraic $K$-theory of $X$ agrees with the homotopy of $\Omega\K(Q\StProj_X)$.
\end{example}

\begin{definition}
    For $\beta\in\Hoh^2(U_{\et},\mu_n)$, let $\T^\beta$ denote the presheaf of spectra
    \begin{equation*}
      V\mapsto\K(\nSetsB_V).
    \end{equation*}
    Define
    \begin{equation*}
        \T_{k}^\beta(V)=\pi_{k}\T^\beta(V),
    \end{equation*}
    and let $\mathcal{T}^{\beta}_k$ be the sheafification of $\T_k^{\beta}$.
\end{definition}

\begin{definition}\label{def:twistedk}
    Similarly, for $\alpha\in\Hoh^2(U_{\et},\Gm)$, let $\K^{\alpha}$ be the presheaf of spectra
    \begin{equation*}
        V\mapsto\K(\StProjA_V),
    \end{equation*}
    with associated homotopy presheaves
    \begin{equation*}
        \K_{k}^{\alpha}(V)=\pi_{k}\K^{\alpha}(V),
    \end{equation*}
    and presheaves $\mathcal{K}^{\alpha}_k$.
\end{definition}

\begin{remark}
    Note that the presheaf of spectra $\K^{\alpha}$ is in some sense the wrong choice of presheaf. The correct version would be to take Thomason-Trobaugh $K$-theory \cite{tt_higher_1990}. However,
    all of the computations in this paper have to do with the \'{e}tale sheafification of $\K^{\alpha}$. Since the two versions agree on affine schemes, it follows that their
    \'{e}tale sheafifications are isomorphic in the homotopy category.
\end{remark}

If $\beta\mapsto\alpha$ in $\Hoh^2(U_{\et},\mu_n)\rightarrow\Hoh^2(U_{\et},\Gm)$, then the twisted unit morphism $i^{\beta}$ of Lemma~\ref{lem:twistedunit} gives a morphism of presheaves of spectra
\begin{equation*}
    \K(i^{\beta}):\T^{\beta}\rightarrow\K^{\alpha}.
\end{equation*}
This map is crucial to the proof of the \textbf{bound} property of the \'etale index.

\begin{lemma}
    Let $\beta\in\Hoh^2(U_{\et},\mu_n)$. Then, the stalk of $\mathcal{T}^{\beta}_j$ at a geometric point $\overline{x}\rightarrow U$ is naturally isomorphic to
    \begin{equation*}
        \pi_j^sB(\mu_n(k(\overline{x})))\oplus\pi_j^s,
    \end{equation*}
    where $k(\overline{x})$ is the (separably closed) residue field of $\overline{x}$, $\pi_j^s$ is the $j$th stable homotopy group of $S^0$, and $BG$ denotes the topological classifying space of a group $G$.
    \begin{proof}
        It is enough to study the stalk $(T^{\beta}_j)_{\overline{x}}$, as this is isomorphic to $(\mathcal{T}^{\beta}_j)_{\overline{x}}$. Since the $K$-theory functor
        preserves filtered colimits, because the classifying space construction does,
        \begin{equation*}
            (T^{\beta}_j)_{\overline{x}}    \iso   \colim_{\overline{x}\in V\rightarrow U}T^{\beta}_j(V)   =  \colim_{\overline{x}\in V\rightarrow U}\K_j(\nSetsB_V) \iso \K_j\left(\colim_{\overline{x}\in V\rightarrow U}\nSetsB_V\right).
        \end{equation*}
        But, $\colim_{\overline{x}\in V\rightarrow U}\nSetsB_V$ is equivalent, by the arguments of \cite[EGA IV 8.5]{ega4.3},
        to the category of finite and faithful $\mu_n(\mathcal{O}^{\sh}_{U,\overline{x}})\iso\mu_n(k(\overline{x}))$-sets. Therefore,
        \begin{equation*}
            (\mathcal{T}^{\beta}_j)_{\overline{x}}\iso\K_j(\nSets_{\overline{x}}),
        \end{equation*}
        where $\nSets_{\overline{x}}$ is the symmetric monoidal category of finite and faithful $\mu_n(k(\overline{x}))$-sets and isomorphisms.
        This category is a groupoid equivalent to
        \begin{equation*}
            \coprod_{j\geq 0} S_j\wr \mu_n(k(\overline{x})),
        \end{equation*}
        where $S_j$ is the symmetric group on $j$ letters, and $S_j\wr\mu_n$ is the wreath product.
        The notation means that the stalk is equivalent to the groupoid with connected components indexed by $j\geq 0$, where the automorphism group of an object in the $j$th component is
        \begin{equation*}
            S_j\wr\mu_n(k(\overline{x})).
        \end{equation*}
        Therefore, by the Barratt-Priddy-Quillen-Kahn theorem (see Thomason \cite[Lemma~2.5]{thomason_first_1982}),
        the $K$-theory spectrum of this symmetric monoidal category is weak equivalent to the suspension spectrum $\Sigma^{\infty}(B\mu_n(k(\overline{x})))_+$
        of the classifying space of $B\mu_n(k(\overline{x}))$ with a disjoint basepoint.
        This spectrum is weakly equivalent to $\Sigma^{\infty}\left(B\mu_n(k(\overline{x}))\vee S^0\right)$. This completes the proof.
    \end{proof}
\end{lemma}

If $n$ is prime to the characteristic of $k(\overline{x})$, then $\mu_n(k(\overline{x}))\iso\ZZ/(n)$. Otherwise, let $m$ be the largest divisor of $n$ that is prime to the characteristic.
Then, $\mu_n(k(\overline{x}))\iso\ZZ/(m)$.

\section{Stable homotopy of classifying spaces}\label{sub:stable}

\begin{proposition}\label{prop:kstemsp}
  Let $0<k<2p-3$. Then, the $p$-primary component $\pi_k^s(p)$ of $\pi_k^s$ is zero. And,
  \begin{equation*}
    \pi_{2p-3}^s(p)=\ZZ/(p).
  \end{equation*}
  \begin{proof}
      This follows from the computation of the image of the $J$-morphism (see \cite[Theorem~1.1.13]{ravenel_complex_1986}) and, for example, \cite[Theorem~1.1.14]{ravenel_complex_1986}.
  \end{proof}
\end{proposition}

I thank Peter Bousfield for telling me about the next proposition.
\begin{proposition}\label{nstablehomotopy}
  For $0<k<2p-2$, the stable homotopy group $\pi_k^s(B\ZZ/(p^n))$ is isomorphic to $\ZZ/(p^n)$ for $k$ odd and $0$ for $k$ even.
  \begin{proof}
       Let $p$ be a prime. Recall the stable splitting of Holzsager~\cite{holzsager_stable_1972}
      \begin{equation*}
          \Sigma B\ZZ/(p^n)\riso X_1\vee\cdots\vee X_{p-1},
      \end{equation*}
      where, if $k>0$, the reduced homology of $X_m$ is
      \begin{equation*}
          \tilde{\Hoh}_k(X_m,\ZZ)\riso\begin{cases}
              \ZZ/(p^n) & \text{if $k\cong 2m \mod 2p-2$,}\\
              0         & \text{otherwise.}
          \end{cases}
      \end{equation*}
      Define $C_m$ as the cofiber of
      \begin{equation*}
          M_1\rightarrow X_m,
      \end{equation*}
      where $M_1=M(\ZZ/(p^n),2m)$ is the Moore space with
      \begin{equation*}
          \tilde{\Hoh}_k(M_1,\ZZ)\riso\begin{cases}
              \ZZ/(p^n) & \text{if $k=2m$,}\\
              0         & \text{otherwise,}
          \end{cases}
      \end{equation*}
      when $k>0$.

      The homology of $C_m$ is
      \begin{equation*}
          \tilde{\Hoh}_k(C_m,\ZZ)\riso\begin{cases}
              \ZZ/(p^n) & \text{if $k>2m$ and $k\cong 2m \mod 2p-2$,}\\
              0         & \text{otherwise.}
          \end{cases}
      \end{equation*}
      Therefore, the map
      \begin{equation*}
          M_2=M(\ZZ/(p^n),2m+2p-2)\rightarrow C_m
      \end{equation*}
      is a $(2m+4p-5)$-equivalence. Thus, for $k<2m+4p-5$ (resp. $k=2m+4p-5$), the map
      \begin{equation*}
          \pi_k^s(M_2)\rightarrow\pi_k^s(C_m)
      \end{equation*}
      is an isomorphism (resp. surjection). Therefore, there is an exact sequence
      \begin{equation}\label{eq:exactmoore}
          \begin{split}
              \pi_{2m+4p-5}^s(M_2)  &\rightarrow\pi_{2m+4p-6}^s(M_1)\rightarrow\pi_{2m+4p-6}^s(X_m)\rightarrow\pi_{2m+4p-6}^s(M_2)\rightarrow\cdots\\
                                    &\rightarrow\pi_k^s(M_1)\rightarrow\pi_k^s(X_m)\rightarrow\pi_k^s(M_2) \rightarrow\cdots
          \end{split}
      \end{equation}

      Let $M(\ZZ/(p^n))$ be the Moore spectrum. It is the cofiber of the multiplication by $p^n$ map on the sphere spectrum $S$. Thus, its stable homotopy groups fit into exact sequences
      \begin{equation*}
          0\rightarrow \pi_k^s\otimes_{\ZZ}\ZZ/(p^n)\rightarrow\pi_k(M(\ZZ/(p^n)))\rightarrow\Tor_1^{\ZZ}(\pi_{k-1}^s,\ZZ/(p^n))\rightarrow 0.
      \end{equation*}
      These sequences are in fact split when $p$ is odd or when $p=2$ and $n>1$. The Moore spaces $M_1$ and $M_2$ are the level $2m$ and $(2m+2p-2)$ spaces of $M(\ZZ/(p^n))$. Thus,
      \begin{align*}
          \pi_k^s(M_1)  &=\pi_{k-2m}(M(\ZZ/(p^n)))\\
          \pi_k^s(M_2)  &=\pi_{k-2m-2p+2}(M(\ZZ/(p^n))).
      \end{align*}

      By Proposition~\ref{prop:kstemsp}, the first $p$-torsion in $\pi_k^s$ is a copy of $\ZZ/(p)$ in degree $k=2p-3$. Therefore, the first two non-zero stable homotopy groups of $M_1$ and $M_2$ are
      \begin{align*}
          \pi_{2m}^s(M_1)       &=\ZZ/(p^n)\\
          \pi_{2m+2p-3}^s(M_1)  &=\ZZ/(p)\\
          \pi_{2m+2p-2}^s(M_2)  &=\ZZ/(p^n)\\
          \pi_{2m+4p-5}^s(M_2)  &=\ZZ/(p).
      \end{align*}
      Using the exact sequence~\eqref{eq:exactmoore}, it follows that the first non-zero stable homotopy group of $X_m$ is
      \begin{equation*}
          \pi_{2m}^s(X_m)=\ZZ/(p^n).
      \end{equation*}
      The next potentially non-zero stable homotopy group fits into the exact sequence~\eqref{eq:exactmoore} at degree $2m+2p-3$:
      \begin{equation*}
          \ZZ/(p^n)\rightarrow\ZZ/(p)\rightarrow\pi_{2m+2p-3}^s(X_m)\rightarrow 0.
      \end{equation*}
      It follows that
      \begin{equation*}
          \pi_k^s(\Sigma B\ZZ/(p^n))=\begin{cases}
              \ZZ/(p^n) &   \text{if $0<k<2p-1$ and $k$ is even,}\\
              0         &   \text{if $0<k<2p-1$ and $k$ is odd}.
          \end{cases}
      \end{equation*}
      The theorem follows immediately.
  \end{proof}
\end{proposition}

\begin{corollary}\label{cor:nstablehomotopy}
    If,
    \begin{equation*}
        \ZZ/(n)=\bigoplus_{q|n}\ZZ/(q^{v_q(n)}),
    \end{equation*}
    where $q$ ranges over the prime divisors of $n$, then, for $0<k<2\min_{q|n}(q)-2$, $\pi_k^s(\ZZ/(n))\iso \ZZ/(n)$ when $k$ is odd and $\pi_k^s(B\ZZ/(n))=0$ when $k$ is even.
    \begin{proof}
        This follows from the proposition, since
        \begin{equation*}
        BG\rwe \vee_{q|n} B\ZZ/(q^{v_q(n)}).
        \end{equation*}
    \end{proof}
\end{corollary}

\begin{corollary}\label{cor:stablesheaves}
    Denote by $m_j$ the exponent of the finite abelian group $\pi_j^s$ for $j\geq 1$. If $\beta\in\Hoh^2(U_{\et},\mu_n)$, then, for
    \begin{equation*}
        0<j<2\min_{q|n}(q)-2,
    \end{equation*}
    the cohomology group $\Hoh^k(U_{\et},\mathcal{T}_j^{\beta})$ is annihilated by $n\cdot m_{j}$ when $j$ is odd and by $m_{j}$ when $j$ is even.
    \begin{proof}
        The stalk of $\mathcal{T}_j^{\beta}$ at $\overline{x}\rightarrow U$ is isomorphic to
        \begin{equation*}
            \pi_{j}^s(B\mu_n(k(\overline{x})))\oplus\pi_{j}^s.
        \end{equation*}
        But, $\mu_n(k(\overline{x}))\iso\ZZ/(m)$, where $m$ is the largest divisor of $n$ prime to the characteristic of $k(\overline{x})$.
        The corollary now follows from the computation of Corollary~\ref{cor:nstablehomotopy}.
    \end{proof}
\end{corollary}

\section{Homotopy sheaves are isomorphic}\label{sub:homotopysheaves}

\begin{proposition}\label{prop:homotopysheaves}
    Fix an element $\alpha\in\Hoh^2(U_{\et},\Gm)$. Then, for all $n\geq 0$, the homotopy sheaves $\mathcal{K}^\alpha$ and $\mathcal{K}$ are naturally isomorphic. Similarly, if $\beta\in\Hoh^2(U_{\et},\mu_n)$,
    then $\mathcal{T}^\beta\iso\mathcal{T}$.
    \begin{proof}
        Here is a proof for the case of $\alpha\in\Hoh^2(U_{\et},\Gm)$. The proof of the other case is identical.
        
        Let $\mathcal{U}=(U_i)_{i\in I}\rightarrow U$ be a cover over which $\alpha$ is trivial (this is possible by the local triviality of sheaf cohomology). 
        Then, by Lemma~\ref{lem:linebundles}, there are $\alpha$-twisted line bundles $\mathcal{L}_i$ on each $U_i$.
        These define equivalences of stacks $\theta_i:\StProj|_{U_i}\rightarrow\StProjA|_{U_i}$ for all $i$ given by
        \begin{equation*}
          \theta_i(V)(\mathcal{P})=\mathcal{L}_i\otimes\mathcal{P},
        \end{equation*}
        when $V\rightarrow U_i$.
        These equivalences induce point-wise weak equivalences of $K$-theory presheaves:
        $\theta_i:\K|_{U_i}\rightarrow\K^{\alpha}|_{U_i}$. This means that for all \'etale maps $V\rightarrow U_i$,
        \begin{equation*}
            (\theta_i)|_{V}:\K|_V\rightarrow\K^{\alpha}|_V
        \end{equation*}
        is a weak equivalence. It follows that on $U_i$ there are isomorphisms of homotopy presheaves:
        \begin{equation*}
            \theta_i:(\K_n)|_{U_i}\riso(\K_n^\alpha)|_{U_i}.
        \end{equation*}
        In fact, the $\theta_i$ glue at the level of homotopy sheaves.
        It suffices to check that, on $U_{ij}=U_i\times_{U}U_j$, the auto-equivalence of $\StProj|_{U_{ij}}$ given by tensoring by $\mathcal{M}_{ij}=\mathcal{L}_i^{-1}\otimes\mathcal{L}_j$
        is locally homotopic to the identity. But, there is a trivialization of $\mathcal{M}_{ij}$, over a cover $\mathcal{V}$ of $U_{ij}$. So, on each element $V$ of $\mathcal{V}$, there
        is an isomorphism $\sigma_V:\mathcal{O}_{U_V}\riso\mathcal{M}_{ij}|_{V}$. This induces a natural transformation from the identity to $\theta_i^{-1}\circ\theta_j$ on $V$.
        But, the $\K$-functor takes natural transformations to homotopies of maps of spectra.
        So, on $V$, $\theta_i|_V=\theta_j|_V:(\K_n)|_V\rightarrow(\K_n^\alpha)|_V$.
        It follows that the $\theta_i$ glue to give isomorphisms \emph{of sheaves}
        \begin{equation*}
            \theta:\mathcal{K}_n\riso\mathcal{K}_n^{\alpha},
        \end{equation*}
        as desired.
    \end{proof}
\end{proposition}

\section{The period-index problem}\label{sec:periodindex}

\begin{definition}\label{def:etalesheaf}
    Let $\K^{\alpha,\et}$ (resp. $\T^{\beta,\et}$) denote the \'etale sheafification of $\K^{\alpha}$ (resp. $\T^{\beta}$) with respect to the local model structure on presheaves of spectra.
    This is the model structure in which cofibrations are given by cofibrations of spectra in the sense of Bousfield and Friedlander~\cite{bousfield_homotopy_1978}, and weak equivalences are morphisms
    that induce isomorphisms of all homotopy sheaves. Since $U$ is of finite cohomological dimension,
    specific models are given by Thomason~\cite[Definition~1.33]{thomason_algebraic_1985}. There are convergent spectral sequences, called Brown-Gersten or descent spectral sequences,
    \begin{align}\label{eq:ss}
        \Eoh_2^{s,t}  &=    \Hoh^s(U_{\et},\mathcal{K}_t^{\alpha})  \Rightarrow\K_{t-s}^{\alpha,\et}(U)\\
        \Eoh_2^{s,t}  &=    \Hoh^s(U_{\et},\mathcal{T}_t^{\beta})   \Rightarrow\T_{t-s}^{\beta,\et}(U)
    \end{align}
    with differentials $d_k^{\alpha}$ of degree $(k,k-1)$; see~\cite[Proposition~1.36]{thomason_algebraic_1985}.
\end{definition}

\begin{definition}\label{def:eti}
    Let $\alpha\in\Hoh^2(U_{\et},\Gm)_{\tors}$. Define the \'etale index of $\alpha$, $eti(\alpha)$,
    to be the positive generator of the image of the edge map (or rank map) $\K_0^{\alpha,\et}(U)\rightarrow\Hoh^0(U_{\et},\mathcal{K}_0^{\alpha})\iso\ZZ$ in the descent spectral sequence.
\end{definition}

\begin{remark}
    The map of presheaves $\K_0^{\alpha,\et}\rightarrow\ZZ$ is called the rank map because the composite $\K_0^{\alpha}\rightarrow\K_0^{\alpha,\et}\rightarrow\ZZ$ is the usual rank map
    on the presheaf of $\alpha$-twisted Grothendieck groups.
\end{remark}

\begin{lemma}[\textbf{Computability}]\label{lem:computability}
    Let $\alpha\in\Hoh^2(U_{\et},\Gm)_{\tors}$. Then, $eti(\alpha)$ is the unique smallest positive integer in $\Hoh^0(U_{\et},\mathcal{K}_0^{\alpha})\iso\ZZ$ such that
    \begin{equation*}
      d_k^\alpha(eti(\alpha))=0
    \end{equation*}
    for all $k\geq 2$.
    \begin{proof}
        This follows immediately from the convergence of the descent spectral sequence~\eqref{eq:ss}.
    \end{proof}
\end{lemma}

\begin{lemma}[\textbf{Obstruction}]\label{lem:obstruction}
    For $\alpha\in\Hoh^2(U_{\et},\Gm)_{\tors}$,
    \begin{equation*}
        eti(\alpha)|deg(\mathcal{A})
    \end{equation*}
    for any Azumaya algebra $\mathcal{A}$ in the class of $\alpha$.
    \begin{proof}
        Suppose that $\mathcal{A}$ is in the class of $\alpha$ and that $m=deg(\mathcal{A})$. Then, by Lemma~\ref{lem:representation}, there is an $\alpha$-twisted lffr sheaf of rank $m$.
        Hence, $m$ is in the image of $rank:\K_0^{\alpha}(U)\rightarrow\ZZ$. Since the rank homomorphism factors through $\K_0^{\alpha,\et}(U)\rightarrow\ZZ$, the lemma follows from
        the definition of the \'etale index.
    \end{proof}
\end{lemma}

\begin{theorem}[\textbf{Divisibility}~\cite{antieau_cech_2010}]
    For $\alpha\in\Hoh^2(U_{\et},\Gm)_{\tors}$,
    \begin{equation*}
        per(\alpha)|eti(\alpha).
    \end{equation*}
\end{theorem}

\begin{example}
    If $D$ is a cyclic division algebra $(x,y)_{\zeta_n}$ over a field of characteristic prime to $n$, so that $per(D)=ind(D)=n$, then $eti(D)=n$.
\end{example}

\begin{example}
    If $D/k$ is a division algebra, and if $l/k$ is a finite separable field extension of degree prime to $per(D)$, then a standard argument using norm maps says that $eti(D_l)=eti(D)$.
\end{example}

\begin{example}
    Let $Q$ be the non-separated quadric with $\alpha$ the non-zero cohomological Brauer class \cite{edidin_brauer_2001}. Then $per(\alpha)=eti(\alpha)=2$, while $ind(\alpha)=+\infty$.
\end{example}

Denote by $m_j$ the exponent of $\pi_j^s$, the $j$th
stable homotopy group of $S^0$, and let $n_j^\alpha$ denote the exponent of $\pi_j^s(B\ZZ/(per(\alpha)))$. Finally,
let $l_j^\alpha$ denote the exponent of $\pi_j^s\oplus\pi_j^s(B\ZZ/(per(\alpha)))$. So, $l_j^\alpha$ is the least common multiple of $m_j$ and $n_j^\alpha$.

\begin{theorem}[\textbf{Bound}]\label{thm:etibound}
  Let $U$ be a connected scheme of cohomological dimension $d$. Let $\alpha\in\Hoh^2(U_{\et},\Gm)_{\tors}$ be such that $per(\alpha)$ is prime to the characteristic of all residue fields of $U$. Then,
  \begin{equation*}
      eti(\alpha)|\prod_{j\in\{1,\ldots,d-1\}} l_j^\alpha.
  \end{equation*}
  \begin{proof}
    Because of the assumption on $per(\alpha)$ and the residue characteristics of $U$, the sequence of sheaves
    \begin{equation*}
        1\rightarrow\mu_{per(\alpha)}\rightarrow\Gm\xrightarrow{per(\alpha)}\Gm\rightarrow 1
    \end{equation*}
    is exact. Thus, there is a lift $\beta$ of $\alpha$ in $\Hoh^2(U_{\et},\mu_{per(\alpha)})$. There is a morphism of descent spectral sequences \cite{thomason_algebraic_1985}
    \begin{equation*}
        \Hoh^s(U_{\et},\mathcal{T}_t^{\beta})\rightarrow\Hoh^s(U_{\et},\mathcal{K}_t^{\alpha})
    \end{equation*}
    induced by $\K(i^{\beta}):\T^{\beta}\rightarrow\K^{\alpha}$. Let $d_k^\beta$ denote the $k$th differential in the descent spectral sequence for $\T^\beta$.
    As the class $1\in\Hoh^0(U_{\et},\mathcal{T}_0^\beta)$ maps to the class $1\in\Hoh^0(U_{\et},\mathcal{K}_0^\alpha)$,
    if $d_k^\beta(m)=0$ for $2\leq k\leq k'$, then $d_k^\alpha(m)=0$ for $2\leq k\leq k'$. The differential $d_k^\beta$ lands in a subquotient of $\Hoh^k(U,\mathcal{T}_{k-1}^\beta))$.
    Therefore, $d_k^\beta$ lands in a group of exponent at most $l_{k-1}^\alpha$, by Corollary~\ref{cor:stablesheaves}.
    Since the sheaves $\mathcal{T}_k^\beta$ are torsion for $k>0$, the differentials $d_k^\beta$ vanish for $k>d$.
  \end{proof}
\end{theorem}

\begin{definition}
    Let $K$ be a field, and let $S$ be a non-empty set of primes. Let $cd_S k$ be the supremum of all the cohomological dimensions $cd_q k$ for all primes $q\in S$.
\end{definition}

\begin{theorem}\label{thm:periodindex}
    Let $K$ be a field, and let $\alpha\in\Br(K)=\Hoh^2(K,\Gm)$ be such that $n=per(\alpha)$ is prime to the characteristic of $K$.
    Let $S$ be the set of prime divisors of $n$, and suppose that $d=cd_S k<2\min_{q\in S}(q)$. Then,
    \begin{equation*}
        eti(\alpha)|(per(\alpha))^{\lfloor\frac{d}{2}\rfloor}.
    \end{equation*}
    \begin{proof}
        Set $c=\lfloor\frac{d}{2}\rfloor$. Combining Theorem~\ref{thm:etibound} and Corollary~\ref{cor:stablesheaves}, it follows that, if $d$ is even, then
        \begin{equation*}
            d_k^{\beta}(an^c)=0
        \end{equation*}
        for all $k\geq 2$, where $a$ is prime to $n$. The same reasoning shows that if $d$ is odd, then
        \begin{equation*}
            d_k^{\beta}(an^c)=0
        \end{equation*}
        when $2\leq k\leq d-1$. By \cite{suslin_local_1984}, the stalks of $\mathcal{K}_{2j}^{\alpha}$ are torsion-free for $j>0$. Therefore, the maps
        \begin{equation*}
            \Hoh^m(K,\mathcal{T}_{2j})\rightarrow\Hoh^m(K,\mathcal{K}_{2j})
        \end{equation*}
        are zero for $j>0$ and all $m$. It follows that if $d_k^{\beta}(m)=0$ for $2\leq k\leq 2j$, then $d_k^{\alpha}(m)=0$ for $2\leq k\leq 2j+1$. Therefore, when $d$ is odd,
        \begin{equation*}
            d_k^{\alpha}(an^c)=0
        \end{equation*}
        for $2\leq k\leq d$ and hence for all $k\geq 2$.

        Thus,
        \begin{equation*}
            eti(\alpha)|an^c,
        \end{equation*}
        where $a$ is relatively prime to $n$.
        On the other hand, as $K$ is a field, the primes divisors of $per(\alpha)$ and $eti(\alpha)$ are the same since $eti(\alpha)|ind(\alpha)$. So,
        \begin{equation*}
        eti(\alpha)|n^f
        \end{equation*}
        for some positive integer $f$. It follows that
        \begin{equation*}
            eti(\alpha)|n^{\min(c,f)}|n^c.
        \end{equation*}
        This completes the proof.
   \end{proof}
\end{theorem}

The condition $d<2\min_{S}(q)$ excludes no primes for function fields of curves, surfaces, or three-folds. It excludes the prime $2$ for function fields of four-folds and five-folds.

The \textbf{bound} property and the method of the proof of Theorem~\ref{thm:periodindex} can be used to give bounds on $eti(\alpha)$ whenever the stable homotopy is known in a sufficiently
large range. But, the exponent $\lfloor\frac{d}{2}\rfloor$ will no longer suffice (with this method). For instance, if $k$ is such that $cd_2 k=4$ and $k$ is not characteristic $2$,
then for any $\alpha\in\Br(k)$ of $per(\alpha)=2$, these arguments give $eti(\alpha)|per(\alpha)^4$. The extra factor of $per(\alpha)^2$ comes from the fact that $\pi_3^s=\ZZ/(24)$.

Let
\begin{align*}
    \K_0^{\alpha}(X)^{(0)}  &=  \K_0^{\alpha}/\ker\left(\K_0^{\alpha}(X)\xrightarrow{rank}\ZZ\right)\\
    \K_0^{\alpha,\et}(X)^{(0)}  &=  \K_0^{\alpha,\et}/\ker\left(\K_0^{\alpha,\et}(X)\xrightarrow{rank}\ZZ\right).
\end{align*}
When $\alpha$ is trivial, the natural inclusion
\begin{equation}\label{eq:weight0}
    \K_0^{\alpha}(X)^{(0)}\rightarrow\K_0^{\alpha,\et}(X)^{(0)}
\end{equation}
is an isomorphism.

\begin{corollary}
    The map of Equation~\eqref{eq:weight0} is not surjective in general when $\alpha$ is not trivial.
    \begin{proof}
        For example, let $k(C)$ be the function field of a curve over a $p$-adic field. Jacob and Tignol have shown in an appendix of \cite{saltman_division_1997}
        that there are division algebras over $k(C)$ for which
        $ind(\alpha)=per(\alpha)^2$. However, since these fields are of cohomological dimension $3$, it follows that $eti(\alpha)=per(\alpha)$. Thus, the map is not surjective for
        $X=\Spec k(C)$.
    \end{proof}
\end{corollary}

\begin{conjecture}\label{conj:etii}
    Let $k=\CC((t_1))\cdots((t_d))$ be an iterated Laurent series field over the complex numbers. Then, for $\alpha\in\Br(k)$,
    \begin{equation*}
        eti(\alpha)=ind(\alpha).
    \end{equation*}
\end{conjecture}

One reason to believe this conjecture is that for $d$-local fields $k$, Becher and Hoffman have established \cite{becher_symbol_2004} that the index satisfies
\begin{equation*}
    ind(\alpha)|per(\alpha)^{\lfloor\frac{d}{2}\rfloor},
\end{equation*}
for all $\alpha\in\Br(k)$.

\bibliographystyle{amsplain-pdflatex}
\bibliography{brauer,mypapers}

\noindent
Benjamin Antieau
[\texttt{\href{mailto:antieau@math.ucla.edu}{antieau@math.ucla.edu}}]\\
UCLA\\
Math Department\\
520 Portola Plaza\\
Los Angeles, CA 90095-1555\\

\end{document}